\documentclass[aos,MSNbibl,nameyear,dvips]{arximspdf}
\usepackage{graphicx}
\doi{10.1214/13-AOS1175D} 
\referstodoi{10.1214/13-AOS1175}
\volume{42}
\issue{2}
\pubyear{2014}
\firstpage{493}
\lastpage{500}

\makeatletter
\newcommand{\Exp}{\operatorname{Exp}}
\newcommand{\supp}{\operatorname{supp}}
\makeatother

\begin{document}
\begin{frontmatter}

\title{Discussion: ``A significance test for the lasso''\thanksref{T1}}
\runtitle{Discussion}

\begin{aug}
\author[a]{\fnms{Jinchi} \snm{Lv}\corref{}\ead[label=e1]{jinchilv@marshall.usc.edu}}
\and
\author[b]{\fnms{Zemin} \snm{Zheng}\ead[label=e2]{zeminzhe@usc.edu}}
\runauthor{J. Lv and Z. Zheng}
\affiliation{University of Southern California}
\address[a]{Data Sciences and Operations Department\\
Marshall School of Business\\
University of Southern California\\
Los Angeles, California 90089\\
USA\\
\printead{e1}} 
\address[b]{Department of Mathematics\\
University of Southern California\\
Los Angeles, California 90089\\
USA\\
\printead{e2}}
\pdftitle{Discussion of ``A significance test for the lasso''}
\end{aug}
\thankstext{T1}{Supported by NSF CAREER Award DMS-09-55316.}

\received{\smonth{12} \syear{2013}}



\end{frontmatter}

Professors Lockhart, Taylor, Tibshirani and Tibshirani are to be
congratulated for their innovative and valuable contribution to the
important and timely problem of testing the significance of covariates
for the Lasso. Since the invention of the Lasso in \citet{tibshirani2}
for variable selection, there has been a huge growing literature
devoted to its theory and implementation, its extensions to various
model settings and different variants and developing more general
regularization methods. Most of existing studies have focused on the
prediction, estimation and variable selection properties ranging from
consistency in prediction and estimation to consistency in model
selection in terms of recovery of the true underlying sparse model. The
problem of deriving the asymptotic distributions for regularized
estimators, as the global or computable solutions, in high dimensions
is relatively less well studied.

How to develop efficient significance testing procedures for the
regularization methods is particularly important since in real
applications one would like to assess the significance of selected
covariates with their $p$-values. Such $p$-values are also crucial for
multiple comparisons in testing the significance of a large number of
covariates simultaneously. In contrast to the use of some resampling or
data splitting techniques for evaluating the significance, in the
present paper Lockhart, Taylor, Tibshirani and Tibshirani propose a
novel powerful yet simple covariance test statistic $T_k$ for testing
the significance of the covariate $X_j$ that enters the model at the
$k$th step of the piecewise linear Lasso solution path in the linear
regression model setting. Such a test statistic is shown to have an
exact $\Exp(1)$ asymptotic null distribution in the case of
orthonormal design matrix and the case of $k = 1$ (i.e., the global
null with zero true regression coefficient vector) for general design
matrix. In the general case, the $\Exp(1)$ distribution provides a
conservative asymptotic null distribution. The significance test for
the Lasso proposed in the paper is elegant thanks to its simplicity and
theoretical guarantees in high dimensions.\looseness=-1

We appreciate the opportunity to comment on several aspects of this
paper. In particular, our discussion will focus on four issues: (1)
alternative test statistics, (2)~the event $B$ and generalized
irrepresentable conditions, (3) model misspecification, and (4) more
general regularization methods.

\section{Alternative test statistics} \label{Sec1}
The covariance test statistic $T_k$ associated with covariate $X_j$ is
defined as the covariance between the response vector $y$ and the net
contribution of covariate $X_j$ toward the mean vector $X \beta$ at
the $(k+1)$th step of the Lasso solution path with regularization
parameter $\lambda= \lambda_{k + 1}$, scaled by the inverse of the
error variance $\sigma^2$. Since the Lasso solution is gauged by the
regularization parameter $\lambda$, a key ingredient in the definition
of $T_k$ is a refitting of the Lasso on the previous support at the
$k$th step right before the inclusion of covariate $X_j$ with the
reduced regularization parameter $\lambda= \lambda_{k + 1}$. This
alignment of the regularization parameter yields a more accurate
account of the contribution of covariate $X_j$ conditional on
previously selected covariates before the next step occurs (either an
addition or a deletion of a covariate).

In view of the geometrical representation of the Lasso solution path,
the choice of a common regularization parameter amounts to that of a
common correlation in magnitude between selected covariates and the
residual vector, provided that all covariate vectors are aligned to a
common scale. In this sense, the covariance test statistic $T_k$ bears
some similarity to the conventional chi-squared test statistic, in
terms of the reduction of the residual sum of squares, for evaluating
the significance of the contribution of a newly added covariate. A main
difference is that the above correlation is constrained as a fixed
number zero in the latter, while it is adaptively
determined at the $(k+1)$th step in the Lasso.

From the estimation point of view, it seems appealing to impose a
smaller regularization to reduce the bias issue incurred by shrinkage.
The bias issue can also affect the significance of relatively weak
covariates. Therefore, a natural extension of the covariance test
statistic $T_k$ can be
%
\begin{equation}
\label{001} T_{k, c} = \bigl(\bigl\langle y, X \hat{\beta}(c
\lambda_{k + 1})\bigr\rangle- \bigl\langle y, X_A \tilde{
\beta}_A(c \lambda_{k + 1})\bigr\rangle\bigr)/
\sigma^2
\end{equation}
for some constant $0 \leq c \leq1$, where $\hat{\beta}(c \lambda_{k
+ 1})$ and $\tilde{\beta}_A(c \lambda_{k + 1})$ are the Lasso
estimators with regularization parameter $c \lambda_{k + 1}$
constrained on sets of covariates $A \cup\{j\}$ and $A$, respectively.
Clearly, $T_{k, c}$ with $c = 1$ reduces to $T_k$. An interesting
question is whether a suitable choice of the constant $c$ may lead to
an exact $\Exp(1)$ asymptotic null distribution for the test statistic
$T_{k, c}$ in the case of general design matrix, as opposed to a
conservative exponential limit for $T_k$.

The significance test with the covariance test statistic $T_k$ is a
sequential procedure which evaluates the significance of any newly
selected covariate in the current Lasso\vadjust{\goodbreak} model. It would also be
appealing to test the significance of each active covariate $X_\ell$
conditional on the set of all remaining active covariates $A_\ell= (A
\cup\{j\}) \setminus\{\ell\}$, since some previously significant
covariates may no longer be significant as new covariates enter the
model. A natural procedure seems to be applying the covariance test
statistic $T_k$ or $T_{k, c}$ to each covariate $X_\ell$ with set $A$
replaced by $A_\ell$. This may also be used to test the significance
of covariates in a general Lasso model given by a prespecified
regularization parameter $\lambda$.

\section{The event $B$ and generalized irrepresentable conditions}\label{Sec2}
There seem to be two key conditions for establishing the conservative
exponential limit for the covariance test statistic $T_k$ in the case
of general design matrix: the event $B$ and generalized irrepresentable
conditions. It would be appealing to verify whether these conditions
hold for a given design matrix. The event $B$ condition assumes that
there exist an integer $k_0 \geq0$ and a fixed set $A_0$ of covariates
containing the set $A^* = \supp(\beta^*)$ of all true active
covariates such that with asymptotic probability one, the Lasso model
$A$ at step $k_0$ of the Lasso solution path is identical to $A_0$. In
other words, this condition requires the sure screening property to
hold for the Lasso.

\begin{figure}

\includegraphics{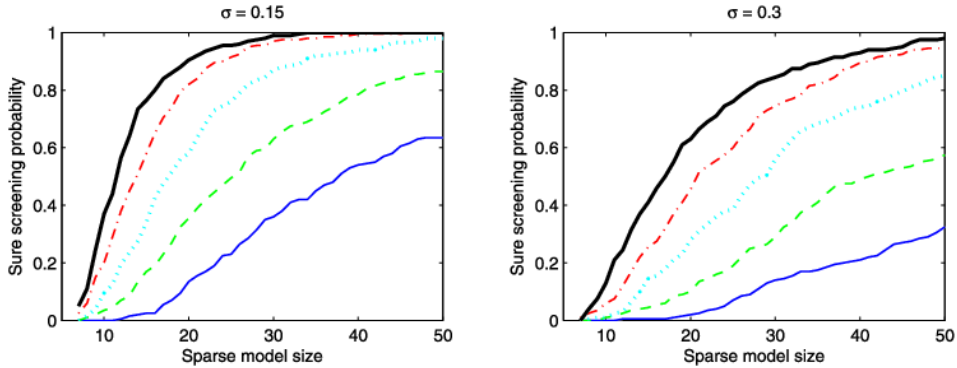}

\caption{Sure screening probability curves as a function of sparse
model size for Lasso with $n = 80$ (thin solid), $90$ (dashed), $100$
(dotted), $110$ (dash--dot) and $120$ (thick solid). Left panel for
$\sigma= 0.15$ and right panel for $\sigma= 0.3$.}
\label{fig1}%
\end{figure}

To provide some insights into the event $B$ condition, let us consider
the setting of linear regression model
%
\begin{equation}
\label{002} y = X \beta^* + \varepsilon
\end{equation}
as in \citet{LF09}. Set $p = 1000$ with true regression
coefficient vector $\beta^* = (1, -0.5, 0.7, -1.2, -0.9, 0.3, 0.55, 0, \ldots, 0)^T$, sample\vspace*{1pt} the rows of the design matrix $X$ as i.i.d.
copies from $N(0, \Sigma)$ with $\Sigma= (0.5^{|i-j|})_{i, j = 1,
\ldots, p}$, and generate error vector $\varepsilon$ independently
from $N(0, \sigma^2 I_n)$ with $\sigma= 0.15$ or $0.3$. Note that the
minimum signal strength is the same as or twice the error standard
deviation. We generated 200 data sets from this model with sample size
$n$ ranging from 80 to 120 and applied the Lasso with the LARS
algorithm [\citet{EHJT04}] to generate the solution path. Figure~\ref{fig1} depicts the sure screening probability curves as a function
of sparse model size for Lasso with different sample size and error
level. The sure screening probability provides an upper bound on the
probability of event $B$. We see that both the signal strength (in
terms of the sample size and error level) and sparse model size are
crucial for the sure screening property of the Lasso estimator. As the
sparse model size and sample size become larger, the Lasso can have
significant sure screening probability. Such a probability can drop as
the noise level increases. It would be interesting to provide some
theoretical understandings on the impacts of these factors on both sure
screening probability and the probability of event $B$ for Lasso.

The generalized irrepresentable condition introduced in the paper
extends the irrepresentable condition [\citet{ZY06}] for
characterizing the model selection consistency of the Lasso, which
means that the true underlying sparse model $A^*$\vadjust{\goodbreak} is exactly recovered
with asymptotic probability one. It involves the fixed set $A_0 \supset
A^*$ introduced in the event $B$ condition. Intuitively, this condition
puts some constraint on the correlation between the noise covariates
and true ones. See, for example, \citet{LF09} and \citet
{FL11}, for examples, and more discussions on these types of conditions
for characterizing the model selection consistency of a wide class of
regularization methods including Lasso.

\section{Model misspecification} \label{Sec3}
The event $B$ condition makes an implicit assumption on the minimum
signal strength. It would be interesting to investigate the more
general case of strong and weak covariates, in which some covariates
may have relatively weak contributions to the response. In such a case,
the true underlying sparse model may no long be contained somewhere on
the solution path given by a regularization method. In other words, the
true model may not be included in the sequence of sparse candidate
models, leading to model misspecification. Apart from missing some true
covariates, model misspecification can generally occur when one
misspecifies the family of distributions. Since model misspecification
is unavoidable in practice, it would be helpful to understand its
impact on statistical inference. For example, \citet{LL13}
recently revealed that the covariance contrast matrix between the
covariance structures in the misspecified model and in the true model
plays a pivotal role in characterizing the impact of model
misspecification on the problem of model selection. It would be
interesting to study the effects of model misspecification in the
context of significance testing.\looseness=-1

\section{More general regularization methods} \label{Sec4}
A key ingredient that makes the covariance test statistic $T_k$ admit
the nice $\Exp(1)$ asymptotic null distribution is the shrinkage
effect induced by the $L_1$-penalty in Lasso which offsets the
inflated\vadjust{\goodbreak}
stochastic variability due to the adaptivity in variable selection.
Many other regularization methods, including concave ones such as the
SCAD [\citet{FL01}], MCP [\citet{Zhang10}] and SICA [\citet{LF09}],
have been proposed for variable selection. A natural and important
question is whether a similar test statistic can be constructed for
testing the significance of covariates for the class of concave
regularization methods. Since these methods are generally nonconvex, it
would be crucial to study the regularized estimate as the global or
computable solution. Consider a fixed regularization parameter $\lambda
$ associated with a penalty function $p_\lambda(t)$ defined on $[0,
\infty)$. One possible test statistic is to extend $T_k$ or $T_{k, c}$
by replacing the constrained Lasso estimators with the corresponding
constrained regularized estimators with the same regularization
parameter. An interesting open question is whether such a generalized
covariance test statistic would have a similar asymptotic null
distribution as for Lasso or a different asymptotic limit may appear.

To gain some insights into these questions, let us consider a natural
extension of the Lasso, the combined $L_1$ and concave regularization
method introduced in \citet{FL13} and defined as the following
regularization problem:
%
\begin{equation}
\label{003} \min_{\beta\in\mathbb{R}^p} \bigl\{(2n)^{-1} \|y - X
\beta\|_2^2 + \lambda_0\|\beta
\|_1 + \bigl\|p_\lambda(\beta)\bigr\|_1 \bigr\},
\end{equation}
where $\lambda_0 = \tilde{c} \{(\log p)/n\}^{1/2}$ for some positive
constant $\tilde{c}$, $p_\lambda(\beta)$ is a compact notation
denoting $p_\lambda(|\beta|) = (p_\lambda(|\beta_1|), \ldots,
p_\lambda(|\beta_p|))^T$ with $|\beta| = (|\beta_1|, \ldots,\break
|\beta_p|)^T$, and $p_\lambda(t)$, $t \in[0, \infty)$, is a penalty
function indexed by the regularization parameter $\lambda\geq0$. This
approach combines the strengths of both Lasso and concave methods in
prediction and variable selection and has enhanced stability compared
with using concave methods alone. They proved that under mild
regularity conditions, the global and computable solutions can enjoy
oracle inequalities under various prediction and estimation losses in
parallel to those in \citet{BRT09} established for
Lasso and Dantzig selector [\citet{CandesTao07}], but with improved
sparsity. In particular, the combined regularization method admits an
explicit bound on the false sign rate, which can be asymptotically vanishing.

\begin{figure}

\includegraphics{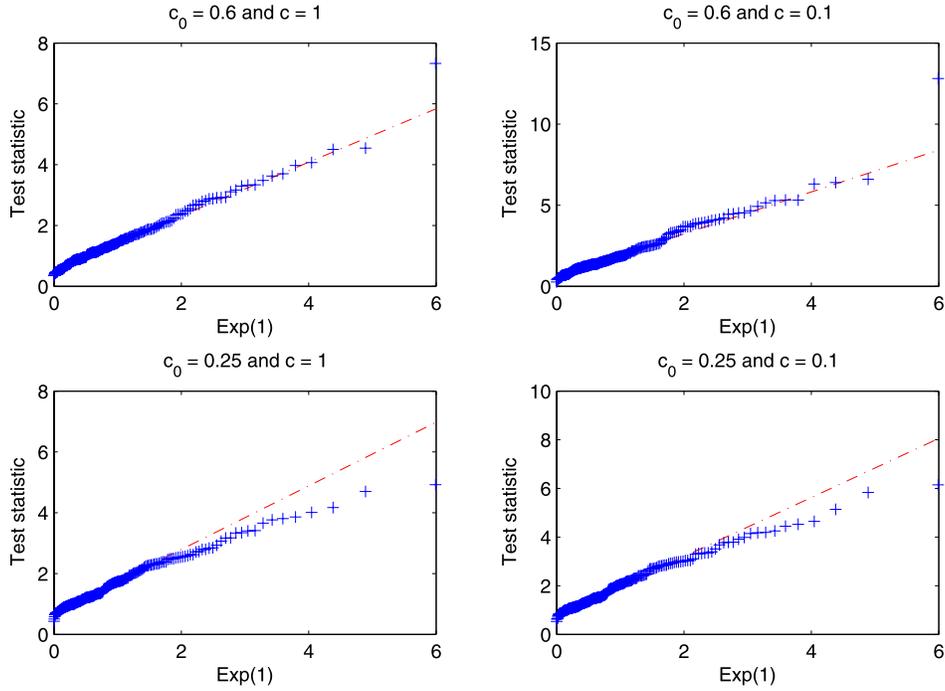}

\caption{Quantile--quantile plots of the covariance test statistic
$T_{k, c}$ with $c = 1$ and $0.1$ versus the $\Exp(1)$ distribution
for the combined $L_1$ and concave regularization with $c_0 = 0.6$ and
$0.25$, respectively.}
\label{fig2}%
\end{figure}
%

\begin{figure}

\includegraphics{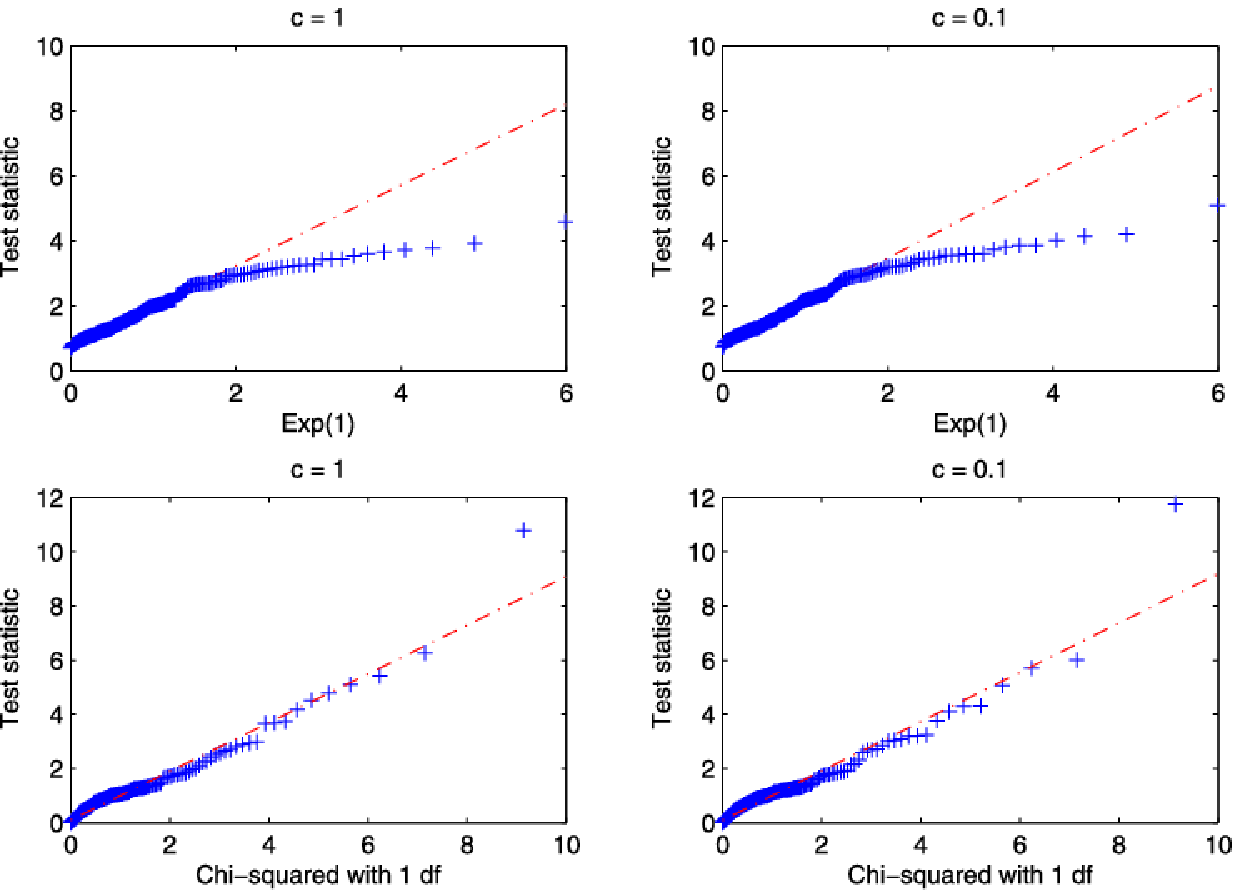}

\caption{Quantile--quantile plots of the covariance test statistic
$T_{k, c}$ with $c = 1$ and $0.1$ for the combined $L_1$ and concave
regularization with $c_0 = 0.1$. Top panel is versus the $\Exp(1)$
distribution and bottom panel is versus the chi-squared distribution
with 1 df.}
\label{fig3}%
\end{figure}

Consider the same example as in Section~\ref{Sec2}, with $n = 100$ and
$\sigma= 0.3$, and use the SICA penalty $p_\lambda(t) = \lambda(a +
1) t/(a + t)$, $t \in[0, \infty)$, with a small shape parameter $a$
for the concave component. For each data set, we applied the combined
$L_1$ and SICA method to generate a sequence of sparse candidate models
with positive constant $\tilde{c} = c_0 \sigma$ and $c_0$ chosen to
be $0.1$, $0.25$ and $0.6$. With tighter control of the false sign
rate, the sure screening property can hold for this method with smaller
sparse model size. Since the true underlying sparse model has seven
variables, we are interested in testing the significance of the eighth
covariate to enter the model. To this end, we used the generalized
covariance test statistics $T_k$ and $T_{k, c}$ (with a slight abuse of
notation) as suggested before, where $T_{k, c}$ with $c = 1$ reduces to
$T_k$ and both $c \lambda_0$ and $c \lambda$ are used in place of
$\lambda_0$ and $\lambda$ in the constrained refittings for~$T_{k,
c}$. The cases of $c = 1$ and $0.1$ were considered.

Figures~\ref{fig2} and \ref{fig3} compare the distributions of the
generalized covariance test statistic $T_{k, c}$ with the $\Exp(1)$
distribution or chi-squared distribution with 1~df over different
combinations of $c_0$ and $c$. Note that as $c_0$ increases, the
combined regularization becomes closer to the Lasso, which is reflected
in Figure~\ref{fig2}. The top panel of Figure~\ref{fig2} for the case
of $c_0 = 0.6$ shows that $\Exp(1)$ fits the distributions of $T_{k,
c}$ well and the bottom panel for the case of $c_0 = 0.25$ suggests
that $\Exp(1)$ is a conservative fit. It is interesting to observe
that the choice of $c = 0.1$ seems to make the distribution of $T_{k,
c}$ closer to the $\Exp(1)$ distribution compared to that of $c = 1$\vadjust{\goodbreak}
in the latter. In the case of $c_0 = 0.1$, the combined regularization
becomes closer to concave regularization. We observe an interesting
transition phenomenon for the distribution of the generalized
covariance test statistic $T_{k, c}$. As demonstrated in Figure~\ref{fig3}, it is now more light-tailed than the $\Exp(1)$ distribution
and interestingly the chi-squared distribution with 1~df provides a
nice fit. It would be interesting to provide theoretical understandings
on such a phenomenon.

\section{Concluding remarks} \label{Sec5}
The covariance test statistic proposed in this paper provides a new
general framework for testing the significance of covariates for the
Lasso and related sparse modeling methods in high dimensions. There are
many interesting and important questions that remain to be answered in
high-dimensional inference. This paper initiates a new area and will
definitely stimulate new ideas and developments in the future. We thank
the authors for their clear and imaginative work.

\section*{Acknowledgment}
We sincerely thank the Co-Editor, Professor Peter Hall, for his kind
invitation to comment on this discussion paper.


%

\printaddresses

\end{document}